\documentclass[12pt,dvips]{amsart}
\usepackage{euler, amsfonts, amssymb, latexsym, epsfig,epic, stmaryrd}

\newtheorem*{Corollary*}{Corollary}
\newtheorem*{Theorem*}{Theorem}

\theoremstyle{remark}

\newtheorem*{Example*}{Example}

\theoremstyle{plain}

\def\<{\langle}
\def\>{\rangle}

\def\nothing{\varnothing}

\newcommand\excise[1]{}

\newcommand{\cellsize}{10}
\newlength{\cellsz} 
\newsavebox{\cell}
\sbox{\cell}{\begin{picture}(\cellsize,\cellsize)
\put(0,0){\line(1,0){\cellsize}}
\put(0,0){\line(0,1){\cellsize}}
\put(\cellsize,0){\line(0,1){\cellsize}}
\put(0,\cellsize){\line(1,0){\cellsize}}
\end{picture}}
\newcommand\cellify[1]{\def\thearg{#1}\def\nothing{}%
\ifx\thearg\nothing
\vrule width0pt height\cellsz depth0pt\else
\hbox to 0pt{\usebox{\cell} \hss}\fi%
\vbox to \cellsz{
\vss
\hbox to \cellsz{\hss$#1$\hss}
\vss}}
\newcommand\tableau[1]{\vtop{\let\\\cr
\baselineskip -16000pt \lineskiplimit 16000pt \lineskip 0pt
\ialign{&\cellify{##}\cr#1\crcr}}}

\begin{document}
\pagestyle{plain}
\mbox{}
\vspace{2ex}

\title{What is a Young Tableau?}
\author{Alexander Yong} 
\address{Department of Mathematics, University of Minnesota,
  Minneapolis, MN 55455, USA; and
  \newline\indent The Fields Institute, 222 College Street, Toronto,
  Ontario, M5T 3J1, Canada}
\email{ayong@math.umn.edu, ayong@fields.utoronto.ca} 
\date{October 30, 2006}
\thanks{I thank Sergey Fomin, Victor Reiner, Hugh Thomas, Alexander Woo
and the editors for helpful suggestions. The author was partially 
supported by NSF grant DMS 0601010 and an NSERC postdoctoral fellowship.}

\maketitle

Young tableaux are ubiquitous combinatorial objects making 
important and inspiring appearances in representation theory, geometry and
algebra. They naturally arise in the study
of symmetric functions, representation theory of the symmetric
and complex general linear groups, and
Schubert calculus of Grassmannians. Discovering and interpreting 
enumerative formulas for Young tableaux (and their generalizations) 
is a core theme of \emph{algebraic combinatorics}. 

Let $\lambda=(\lambda_1\geq \lambda_2\geq \ldots 
\geq \lambda_k\geq 0)$ be a {\bf partition} of size 
$|\lambda|=\lambda_1+\ldots+\lambda_k$, identified
with its {\bf Young diagram}: a left justified shape of $k$ rows
of boxes of length $\lambda_1,\ldots, \lambda_k$. For example,
$\lambda=(4,2,1)$ is drawn
$\tableau{{\ }&{\ }&{\ }&{\ }\\{\ }&{\ }\\{\ }}$. A {\bf (Young) filling} 
of $\lambda$ assigns a positive integer to each box of
$\lambda$, e.g., $\tableau{{2 }&{1 }&{1 }&{4 }\\{6 }&{2 }\\{4 }}$. A filling
is {\bf semistandard} if the entries weakly increase along rows
and strictly increase along columns. A semistandard filling is 
{\bf standard} if it is a bijective assignment of $\{1,2\ldots,|\lambda|\}$.
So $\tableau{{1 }&{2 }&{2 }&{4 }\\{2 }&{3 }\\{4 }}$ is a {\bf semistandard
Young tableau} while $\tableau{{1}&{3}&{4}&{6}\\{2}&{7}\\{5}}$ is a 
{\bf standard Young tableau}, both of {\bf shape} $\lambda$.

We focus on the enumeration and generating series of Young tableaux. 
Frame-Robinson-Thrall's elegant (and nontrivial) {\bf hook-length formula} states that the number of standard Young tableaux of shape $\lambda$ is 
$f^{\lambda}:=\frac{|\lambda|!}{\prod_{b} h_b}$, where the 
product in the denominator is over all boxes $b$ of $\lambda$ and $h_b$ is 
the {\bf hook-length} of $b$, i.e., the number of boxes directly to the right 
or below $b$ (including $b$ itself). 
Thus, 
$f^{(4,2,1)}=\frac{7!}{6\cdot 4 \cdot 2\cdot 1\cdot 3\cdot 1\cdot 1}=35$. 

A similar \emph{hook-content formula} counts the number of semistandard Young 
tableaux, but we now consider instead their 
generating series: fix $\lambda$ and a bound $N$ on the size of
the entries in each semistandard tableau $T$. Let 
${\bf x}^{T}=\prod_{i=1}^{N}x_{i}^{\#{\mbox{\small $i$'s in $T$}}}$.
The {\bf Schur polynomial} is the generating series
$s_{\lambda}(x_1,\ldots, x_N):=\sum_{{\mbox{\small semistandard} \ } T}{\bf x}^{T}$. For example, when $N=3$ and $\lambda=(2,1)$ there are eight semistandard
Young tableaux:
\[\tableau{{1}&{1}\\{2}},\ \  \tableau{{1}&{2}\\{2}},\ \  
\tableau{{1}&{3}\\{2}}, \  \ 
\tableau{{1}&{2}\\{3}}, \ \ 
\tableau{{1}&{1}\\{3}}, \ \ 
\tableau{{1}&{3}\\{3}}, 
\ \ \tableau{{2}&{2}\\{3}},\ \  \tableau{{2}&{3}\\{3}}.\]
The corresponding Schur polynomial, with terms in the same order, is
$s_{(2,1)}(x_1,x_2,x_3)=x_{1}^2 x_2 + x_1 x_2^2 + x_1 x_2 x_3 + x_1 x_2 x_3 + x_1^2 x_3 + x_1 x_3^2
+x_2^2 x_3 + x_2 x_3^2$. In general, these are
{\bf symmetric polynomials}, i.e., 
$s_{\lambda}(x_1,\ldots, x_N)=s_{\lambda}(x_{\sigma(N)},\ldots, x_{\sigma(N)})$
for all $\sigma$ in the symmetric group ${\mathfrak S}_N$ 
(the proof is a ``clever trick'' known
as the \emph{Bender-Knuth involution}).

Both the irreducible complex representations of ${\mathfrak S}_n$
and the irreducible degree $n$ polynomial representations of the general
linear group 
$GL_{N}({\mathbb C})$ are indexed by partitions $\lambda$ with
 $|\lambda|=n$. The associated irreducible ${\mathfrak S}_n$-representation
has dimension equal to $f^{\lambda}$, while the irreducible
$GL_{N}({\mathbb C})$-representation has character 
$s_{\lambda}(x_1,\ldots, x_N)$.
These facts can be proved with an explicit construction of the
respective representations having a basis indexed by the appropriate tableaux.

In algebraic geometry, the Schubert varieties, in the complex Grassmannian 
manifold $Gr(k,{\mathbb C}^n)$ 
of $k$-planes in ${\mathbb C}^n$, are indexed by partitions $\lambda$ contained
inside a $k\times (n-k)$ rectangle. Here, the Schur polynomial 
$s_{\lambda}(x_1,\ldots, x_k)$ represents
the class of the Schubert variety under a natural presentation of the 
cohomology ring $H^{\star}(Gr(k,{\mathbb C}^n))$.

Schur polynomials form a vector space basis (say, over ${\mathbb Q}$) 
of the ring of
symmetric polynomials in the variables $x_1,\ldots, x_N$. Since
a product of symmetric polynomials is symmetric, we can expand the
result in terms of Schur polynomials. 
In particular, define the {\bf Littlewood-Richardson coefficients} $C_{\lambda,\mu}^{\nu}$ by
\begin{equation}
\label{eqn:LR}
s_{\lambda}(x_1,\ldots,x_N)\cdot s_{\mu}(x_1,\ldots, x_N)=
\sum_{\nu}C_{\lambda,\mu}^{\nu} \ s_{\nu}(x_1,\ldots, x_N).
\end{equation}
In fact, $C_{\lambda,\mu}^{\nu}\in{\mathbb Z}_{\geq 0}$! 
These numbers count tensor product multiplicities 
of 
irreducible representations of $GL_{N}({\mathbb C})$. Alternatively, 
they count \emph{Schubert calculus} 
intersection numbers for a 
triple of Schubert varieties in a Grassmannian. 
However, neither of these descriptions of $C_{\lambda,\mu}^{\nu}$ is really a means 
to calculate the number.

The {\bf Littlewood-Richardson rule} combinatorially manifests
the positivity of the $C_{\lambda,\mu}^{\nu}$. Numerous versions
of this rule exist, exhibiting different 
features of the numbers. Here is a standard 
version: take the Young diagram of $\nu$ and remove
the Young diagram of $\lambda$, where the latter is top left justified in the
former (if $\lambda$ is not contained inside 
$\nu$, then declare $C_{\lambda,\mu}^{\nu}=0$); this {\bf skew-shape} is 
denoted $\nu/\lambda$.  Then $C_{\lambda,\mu}^{\nu}$ counts the number of
semistandard fillings $T$
of shape $\nu/\lambda$ such that (a) as the entries are read along rows 
from right to 
left, and from top to bottom, at every point,
the number of $i$'s appearing always is weakly
less than the number of $i-1$'s, for $i\geq 2$; 
and (b) the total number of $i's$ appearing is $\mu_i$. 
Thus $C_{(2,1),(2,1)}^{(3,2,1)}=2$ is witnessed by
$\tableau{{\ }&{\ }&{1}\\{\ }&{1}\\{2}}$ and 
$\tableau{{\ }&{\ }&{1}\\{\ }&{2}\\{1}}$. 

Extending the $GL_{N}({\mathbb C})$ story, 
a generalized Littlewood-Richardson rule exists for 
all complex semisimple Lie groups, where 
\emph{Littelmann paths} generalize Young tableaux~\cite{Littelmann}.
In contrast, the situation is much less satisfactory in 
the Schubert calculus context, although in recent work
with Hugh Thomas \cite{Thomas-Yong}, we made progress for the \linebreak
\emph{(co)minuscule} 
generalization of Grassmannians.

No discussion of Young tableaux is complete without 
the {\bf Schensted correspondence}. This associates each 
$\sigma \in {\mathfrak S}_n$ bijectively with pairs of standard Young tableaux $(T,U)$ of the same shape $\lambda$, where $|\lambda|=n$. This can be used to
 prove the Littlewood-Richardson rule, but is noteworthy in its own right in geometry and representation theory.
 
Given a permutation (in one line notation), e.g., $\sigma=21453\in S_5$, at 
each step $i$ we add a box into some row of the current 
{\bf insertion tableau} ${\widetilde T}$: initially
insert $\sigma(i)$ into the first row of ${\widetilde T}$. 
If no entries $y$ of that row are larger than $\sigma(i)$, 
place $\sigma(i)$ in a new box at 
the end of the row and place a new box containing $i$ at the same place
in the current {\bf recording tableau} ${\widetilde U}$. 
Otherwise, let $\sigma(i)$ replace the leftmost $y>\sigma(i)$ and insert
$y$ into the second row, and so on. This eventually results in two
tableaux of the same shape; Schensted outputs $(T,U)$ after
$n$ steps. In our example, the steps are
\[(\emptyset,\emptyset), (\tableau{{2}},\tableau{{1}}),
\left(\tableau{{1}\\{2}},\tableau{{1}\\{2}}\right),
\left(\tableau{{1}&{4}\\{2}},\tableau{{1}&{3}\\{2}}\right),
\left(\tableau{{1}&{4}&{5}\\{2}},\tableau{{1}&{3}&{4}\\{2}}\right),
\left(\tableau{{1}&{3}&{5}\\{2}&{4}},\tableau{{1}&{3}&{4}\\{2}&{5}}\right)
=(T,U).
\] 
It is straightforward to prove
well-definedness and bijectivity of this procedure. 
Also, $T$ and $U$
encode interesting information about $\sigma$. For example, it is easy to 
show that $\lambda_1$ equals the length 
of the longest increasing subsequence in $\sigma$ (see e.g., work
of Baik-Deift-Johansson~\cite{Baik} for
connections to random matrix theory). A sample harder fact is 
that if $\sigma$ corresponds to $(T,U)$ then $\sigma^{-1}$ corresponds 
to $(U,T)$.

An excellent source for more on the combinatorics of
Young tableaux is~\cite{Stanley}, whereas applications to 
geometry and representation theory are developed 
in~\cite{Fulton}. For a survey containing examples of Young tableaux
for other Lie groups, see~\cite{Sagan}. 
Active research on the topic of Young tableaux
continues, for example, recently in collaboration
with Allen Knutson and Ezra Miller~\cite{Knutson.Miller.Yong}, 
we found 
a simplicial ball of semistandard tableaux, 
together with applications to Hilbert series formulae of 
determinantal ideals.

\end{document}